\documentclass[12pt]{article}
\usepackage{latexsym,amsfonts,amsmath,amssymb}
\newtheorem{theorem}{Theorem}
\newtheorem{corollary}[theorem]{Corollary}
\newtheorem{sublemma}{Lemma}[theorem]
\newtheorem{lemma}[theorem]{Lemma}
\newtheorem{question}[theorem]{Question}
\newtheorem{observation}[theorem]{Observation}
\newtheorem{claim}[theorem]{Claim}
\newtheorem{subclaim}{Claim}[sublemma]
\newtheorem{conjecture}[theorem]{Conjecture}
\newtheorem{fact}[theorem]{Fact}
\newtheorem{definition}[theorem]{Definition}
\newtheorem{remark}[theorem]{Remark}
\newtheorem{example}[theorem]{Example}
\newtheorem{exercise}{Exercise}[section]
\def\Theorem #1.#2 #3\par{\setbox1=\hbox{#1}\ifdim\wd1=0pt
   \begin{theorem}{\rm #2} #3\end{theorem}\else
   \newtheorem{#1}[theorem]{#1}\begin{#1}\label{#1}{\rm #2} #3\end{#1}\fi}
\def\Corollary #1.#2 #3\par{\setbox1=\hbox{#1}\ifdim\wd1=0pt
   \begin{corollary}{\rm #2} #3\end{corollary}\else
   \newtheorem{#1}[theorem]{#1}\begin{#1}\label{#1}{\rm #2} #3\end{#1}\fi}
\def\Lemma #1.#2 #3\par{\setbox1=\hbox{#1}\ifdim\wd1=0pt
   \begin{lemma}{\rm #2} #3\end{lemma}\else
   \newtheorem{#1}[theorem]{#1}\begin{#1}\label{#1}{\rm #2} #3\end{#1}\fi}
\def\SubLemma #1.#2 #3\par{\setbox1=\hbox{#1}\ifdim\wd1=0pt
   \begin{sublemma}{\rm #2} #3\end{sublemma}\else
   \newtheorem{#1}{#1}[theorem]\begin{#1}\label{#1}{\rm #2} #3\end{#1}\fi}
\def\Question #1.#2 #3\par{\setbox1=\hbox{#1}\ifdim\wd1=0pt
   \begin{question}{\rm #2} #3\end{question}\else
   \newtheorem{#1}[theorem]{#1}\begin{#1}\label{#1}{\rm #2} #3\end{#1}\fi}
\def\Observation #1.#2 #3\par{\setbox1=\hbox{#1}\ifdim\wd1=0pt
   \begin{observation}{\rm #2} #3\end{observation}\else
   \newtheorem{#1}[theorem]{#1}\begin{#1}\label{#1}{\rm #2} #3\end{#1}\fi}
\def\Claim #1.#2 #3\par{\setbox1=\hbox{#1}\ifdim\wd1=0pt
   \begin{claim}{\rm #2} #3\end{claim}\else
   \newtheorem{#1}[theorem]{#1}\begin{#1}\label{#1}{\rm #2} #3\end{#1}\fi}
\def\SubClaim #1.#2 #3\par{\setbox1=\hbox{#1}\ifdim\wd1=0pt
   \begin{subclaim}{\rm #2} #3\end{subclaim}\else
   \newtheorem{#1}{#1}[sublemma]\begin{#1}\label{#1}{\rm #2} #3\end{#1}\fi}
\def\Conjecture #1.#2 #3\par{\setbox1=\hbox{#1}\ifdim\wd1=0pt
   \begin{conjecture}{\rm #2} #3\end{conjecture}\else
   \newtheorem{#1}[theorem]{#1}\begin{#1}\label{#1}{\rm #2} #3\end{#1}\fi}
\def\Fact #1.#2 #3\par{\setbox1=\hbox{#1}\ifdim\wd1=0pt
   \begin{fact}{\rm #2} #3\end{fact}\else
   \newtheorem{#1}[theorem]{#1}\begin{#1}\label{#1}{\rm #2} #3\end{#1}\fi}
\def\Definition #1.#2 #3\par{\setbox1=\hbox{#1}\ifdim\wd1=0pt
   \begin{definition}{\rm #2} {\rm #3}\end{definition}\else
   \newtheorem{#1}[theorem]{#1}\begin{#1}\label{#1}{\rm #2} {\rm #3}\end{#1}\fi}
\def\Remark #1.#2 #3\par{\setbox1=\hbox{#1}\ifdim\wd1=0pt
   \begin{remark}{\rm #2} {\rm #3}\end{remark}\else
   \newtheorem{#1}[theorem]{#1}\begin{#1}\label{#1}{\rm #2} {\rm #3}\end{#1}\fi}
\def\Example #1.#2 #3\par{\setbox1=\hbox{#1}\ifdim\wd1=0pt
   \begin{example}{\rm #2} #3\end{example}\else
   \newtheorem{#1}[theorem]{#1}\begin{#1}\label{#1}{\rm #2} #3\end{#1}\fi}
\def\Exercise #1.#2 #3\par{\setbox1=\hbox{#1}\ifdim\wd1=0pt
   {\footnotesize\begin{exercise}{\rm #2} {\rm #3}\end{exercise}}\else
   \newtheorem{#1}[section]{#1}{\footnotesize\begin{#1}\label{#1}{\rm #2} {\rm #3}\end{#1}}\fi}
\def\QuietTheorem #1.#2 #3\par{\setbox1=\hbox{#1}\ifdim\wd1=0pt\proclaim{Theorem {\rm #2}}{#3}\else\proclaim{#1 {\rm #2}}{#3}\fi}
\newcommand{\proclaim}[2]{\smallskip\noindent{\bf #1} {\sl#2}\par\smallskip}
\def\Proclaim #1.#2 #3\par{\proclaim{#1 {\rm #2}}{#3}}
\newenvironment{proof}{\noindent}{\kern2pt\QEDbox\par\bigskip}
\def\Proof#1: {\setbox1=\hbox{#1}\ifdim\wd1=0pt\begin{proof}{\bf Proof: }\else\medskip\begin{proof}{\bf #1: }\fi}
\newcommand{\QED}{\end{proof}}
\def\BF#1.{{\bf #1.}}
\def\Abstract #1\par{\begin{quotation}{\singlespaced\footnotesize{\noindent{\bf Abstract.~}#1}}\end{quotation}}
\def\Title #1\par{\title{#1}\maketitle}
\def\Author #1\par{\author{#1}}
\def\Acknowledgement#1\par{\thanks{#1}}
\def\Chapter #1\par{\chapter{#1}}
\def\Section #1\par{\section{#1}}
\def\QuietSection #1\par{\section*{#1}}
\def\SubSection #1\par{\subsection{#1}}
\def\SubSubSection #1\par{\subsubsection{#1}}
\def\MidTitle #1\par{\bigskip\goodbreak\centerline{\small\bf #1}\bigskip\noindent}
\def\Margin #1\par{\marginpar{\scriptsize #1}}

\newcommand{\singlespaced}{\baselineskip=15pt}
\def\bottomnote #1\par{{\renewcommand{\thefootnote}{}\footnotetext{#1}}}
\def\R{{\mathbb R}}

\newfont{\msam}{msam10 at 12pt}

\newcommand{\of}{\subseteq}
\newcommand{\ofnoteq}{\subsetneq}

\newcommand{\set}[1]{\{\,{#1}\,\}}

\newcommand{\restrict}{\upharpoonright}
\newcommand{\satisfies}{\models}

\newcommand{\converges}{\downarrow}

\newcommand{\union}{\cup}
\newcommand{\Union}{\bigcup}
\newcommand{\intersect}{\cap}

\newcommand{\trianglelt}{\lhd}

\newcommand{\smalllt}{\mathrel{\mathchoice{\raise2pt\hbox{$\scriptstyle<$}}{\raise1pt\hbox{$\scriptstyle<$}}{\scriptscriptstyle<}{\scriptscriptstyle<}}}
\newcommand{\smallleq}{\mathrel{\mathchoice{\raise2pt\hbox{$\scriptstyle\leq$}}{\raise1pt\hbox{$\scriptstyle\leq$}}{\scriptscriptstyle\leq}{\scriptscriptstyle\leq}}}

\newcommand{\leqalpha}{{{\smallleq}\alpha}}

\newcommand{\leqxi}{{{\smallleq}\xi}}

\newcommand{\leqSigma}{{{\smallleq}\Sigma}}

\newcommand{\gcode}[1]{{}^\ulcorner#1{}^\urcorner}
\newcommand{\UnderTilde}[1]{{\setbox1=\hbox{$#1$}\baselineskip=0pt\vtop{\hbox{$#1$}\hbox to\wd1{\hfil$\sim$\hfil}}}{}}
\newcommand{\Undertilde}[1]{{\setbox1=\hbox{$#1$}\baselineskip=0pt\vtop{\hbox{$#1$}\hbox to\wd1{\hfil$\scriptstyle\sim$\hfil}}}{}}
\newcommand{\undertilde}[1]{{\setbox1=\hbox{$#1$}\baselineskip=0pt\vtop{\hbox{$#1$}\hbox to\wd1{\hfil$\scriptscriptstyle\sim$\hfil}}}{}}
\newcommand{\UnderdTilde}[1]{{\setbox1=\hbox{$#1$}\baselineskip=0pt\vtop{\hbox{$#1$}\hbox to\wd1{\hfil$\approx$\hfil}}}{}}
\newcommand{\Underdtilde}[1]{{\setbox1=\hbox{$#1$}\baselineskip=0pt\vtop{\hbox{$#1$}\hbox to\wd1{\hfil\scriptsize$\approx$\hfil}}}{}}

\newcommand{\st}{\mid}
\renewcommand{\th}{{\hbox{\scriptsize th}}}

\def\<#1>{\langle#1\rangle}

\newcommand{\QEDbox}{\fbox{}}

\newcommand{\ORD}{\mathop{\hbox{\sc ord}}}
\newcommand{\WO}{\mathop{\hbox{\sc wo}}}

\newcommand{\NP}{\mathop{\hbox{\it NP}}\nolimits}
\newcommand{\coNP}{\mathop{\hbox{\rm co-\!\it NP}}\nolimits}

\newcommand{\omegaCK}{{\omega_1^{\rm ck}}}
\newcommand{\cell}[1]{\boxit{\hbox to 17pt{\strut\hfil$#1$\hfil}}}
\newcommand{\head}[2]{\lower2pt\vbox{\hbox{\strut\footnotesize\it\hskip3pt#2}\boxit{\cell#1}}}
\newcommand{\boxit}[1]{\setbox4=\hbox{\kern2pt#1\kern2pt}\hbox{\vrule\vbox{\hrule\kern2pt\box4\kern2pt\hrule}\vrule}}
\newcommand{\Col}[3]{\hbox{\vbox{\baselineskip=0pt\parskip=0pt\cell#1\cell#2\cell#3}}}
\newcommand{\tapenames}{\raise 5pt\vbox to .7in{\hbox to .8in{\it\hfill input: \strut}\vfill\hbox to
.8in{\it\hfill scratch: \strut}\vfill\hbox to .8in{\it\hfill output: \strut}}}
\newcommand{\Head}[4]{\lower2pt\vbox{\hbox to25pt{\strut\footnotesize\it\hfill#4\hfill}\boxit{\Col#1#2#3}}}
\newcommand{\Dots}{\raise 5pt\vbox to .7in{\hbox{\ $\cdots$\strut}\vfill\hbox{\ $\cdots$\strut}\vfill\hbox{\
$\cdots$\strut}}}
\renewcommand{\dots}{\raise5pt\hbox{\ $\cdots$}}
\newcommand{\factordiagramup}[6]{$$\begin{array}{ccc}
#1&\raise3pt\vbox{\hbox to60pt{\hfill$\scriptstyle
#2$\hfill}\vskip-6pt\hbox{$\vector(4,0){60}$}}&#3\\ \vbox
to30pt{}&\raise22pt\vtop{\hbox{$\vector(4,-3){60}$}\vskip-22pt\hbox
to60pt{\hfill$\scriptstyle #4\qquad$\hfill}}
     &\ \ \lower22pt\hbox{$\vector(0,3){45}$}\ {\scriptstyle #5}\\
\vbox to15pt{}&&#6\\
\end{array}$$}
\newcommand{\factordiagram}[6]{$$\begin{array}{ccc}
#1&&\\ \ \ \raise22pt\hbox{$\vector(0,-3){45}$}\ {\scriptstyle #2}
&\raise22pt\hbox{$\vector(2,-1){90}$}\raise5pt\llap{$\scriptstyle#3$\qquad\quad}&\vbox
to25pt{}\\ #4&\raise3pt\vbox{\hbox to90pt{\hfill$\scriptstyle
#5$\hfill}\vskip-6pt\hbox{$\vector(4,0){90}$}}&#6\\
\end{array}$$}
\newcommand{\df}{\it} 
\begin{document}
\author{Vinay Deolalikar\\
\normalsize\sc Hewlett-Packard Research\\
\\
Joel David Hamkins\thanks{The research of the second author has been supported in part by grants from
Georgia State University, the Research Foundation of CUNY and the National
Science Foundation.}\\
\normalsize\sc Georgia State University \& The City University of New York\\
\\
Ralf-Dieter Schindler\\
\normalsize\sc Institut f\"ur formale Logik, Universit\"at Wien\\
}

\bottomnote Keywords: Infinite time Turing machines, complexity theory, descriptive set theory. Mathematics
Subject Classification 03D15, 68Q15, 68Q17, 03E15.

\Title $P\not=\NP\intersect \coNP$\/ for Infinite Time Turing Machines

\Abstract Extending results of Schindler \cite{Schindler:Pnot=NP} and Hamkins and Welch
\cite{HamkinsWelch2003:Pf=NPf}, we establish in the context of infinite time Turing machines that $P$\/ is
properly contained in $\NP\intersect \coNP$. Furthermore, $\NP\intersect\coNP$ is exactly the class of
hyperarithmetic sets. For the more general classes, we establish that $P^+=\NP^+\intersect
\coNP^+=\NP\intersect\coNP$, though $P^{++}$ is properly contained in $\NP^{++}\intersect \coNP^{++}$.
Within any contiguous block of infinite clockable ordinals, we show that
$P_\alpha\not=\NP_\alpha\intersect\coNP_\alpha$, but if $\beta$ begins a gap in the clockable ordinals,
then $P_\beta=\NP_\beta\intersect\coNP_\beta$. Finally, we establish that $P^f \neq \NP^f
\intersect\coNP^f$ for most functions $f:\R\to\ORD$, although we provide examples where $P^f = \NP^f
\intersect\coNP^f$ and $P^f \neq \NP^f$.

\Section Introduction

In this article, we take up the question of whether $P=\NP\intersect \coNP$\/ for infinite time Turing
machines. The related $P=\NP$\/ problem was first considered in connection with infinite time Turing
machines by Schindler (the third author) in \cite{Schindler:Pnot=NP}, where he proved that $P\not=\NP$\/
and introduced the other natural complexity classes $P^+$, $\NP^+$, $P^{++}$, $\NP^{++}$, $P_\alpha$,
$\NP_\alpha$, $P^f$ and $\NP^f$. He then showed that $P^+ \not= \NP^+$ and posed the corresponding
questions for $P^{++}$ and $P^f$ when $f$ is a suitable function from $\R$ to the ordinals. Hamkins (the
second author) and Welch answered these questions in \cite{HamkinsWelch2003:Pf=NPf} by showing that
$P^{++}\not=\NP^{++}$ and, more generally, that $P^f\not=\NP^f$ for almost every function $f$. Here, we
extend the analysis of all these complexity classes to the analogues of the question of whether
$P=\NP\intersect \coNP$. Unfortunately, there is no uniform answer, as some of the complexity classes
satisfy the equation and some do not, though the general tendency is towards inequality.

We show, in particular, that $P$\/ is properly contained in $\NP\intersect \coNP$. Furthermore,
$\NP\intersect\coNP$ is exactly the class of hyperarithmetic sets. At the next level, we establish
$P^+=\NP^+\intersect \coNP^+=\NP\intersect\coNP$. At a still higher level, once again $P^{++}$ is properly
contained in $\NP^{++}\intersect \coNP^{++}$. Within any contiguous block of infinite clockable ordinals,
we establish $P_\alpha\not=\NP_\alpha\intersect\coNP_\alpha$, but if $\beta$ begins a gap in the clockable
ordinals, then $P_\beta=\NP_\beta\intersect\coNP_\beta$. Finally, for almost all functions $f:\R\to\ORD$,
the class $P^f$ is properly contained in $\NP^f\intersect \coNP^f$, though there are functions for which
$P^f=\NP^f\intersect \coNP^f$, even with $P^f\not=\NP^f$.

Infinite time Turing machines were introduced by Hamkins and Lewis in
\cite{HamkinsLewis2000:InfiniteTimeTM}, and we refer the reader to that article for reference and
background. Let us quickly describe here for convenience how the machines operate.
\begin{figure}[h]
$$\tapenames\Col111\Col100\Head101s\Col111\Col010\Col101\Col001\Col011\Dots$$
\end{figure}
The hardware of an infinite time Turing machine is the same as that of a classical three-tape Turing
machine: a head moves left and right on a semi-infinite paper tape, reading and writing according to the
rigid instructions of a finite program with finitely many states in exactly the classical manner. The
operation of the machines is extended into transfinite ordinal time by defining the configuration of the
machine at the limit ordinal stages. At such a stage, the head is returned to the leftmost cell, the
machine is placed into the special {\it limit} state, and the tape is updated by placing into each cell the
$\limsup$ of the values appearing in that cell before the limit stage. Thus, if the cell values have
stabilized before a limit, then at the limit the cell displays this stabilized value, and otherwise, when
the cell has changed from $0$ to $1$ and back again unboundedly often before the limit, then at the limit
the cell displays a $1$. Having specified the operation of the machines, one obtains for any program $p$
the corresponding infinite time computable function $\varphi_p$, namely, $\varphi_p(x)=y$ when program $p$
on input $x$ halts with output $y$. The natural context for input and output is infinite binary sequences,
that is, Cantor space ${}^\omega 2$, which we refer to as the set of reals and denote by $\R$. A set
$A\of\R$ is infinite time decidable if its characteristic function is decidable. In the context of certain
time-critical complexity classes, we adopt the formalism for deciding sets with two distinct halt states,
{\it accept} and {\it reject}, so that the machines can announce their decisions as quickly as possible,
without needing to position the head for writing on the output tape. For many of the complexity classes,
however, including $P$, $P^+$, $P^{++}$ and $P_\alpha$ for limit ordinals $\alpha$ and their successors,
the additional steps required for writing on the tape pose no difficulty, and one can dispense with this
formalism in favor of the usual characteristic function notion of decidability.

Many of our arguments will rely on elementary results in descriptive set theory, and we refer readers to
\cite{Moschovakis1980:DescriptiveSetTheory}, \cite{Kechris1995:ClassicalDescripiveSetTheory} and
\cite{MansfieldWeitkamp1985:RecursiveAspectsOfDescriptiveSetTheory} for excellent introductions. For
background material on admissible set theory, we refer readers to
\cite{Barwise1975:AdmissibleSetsStructures}. We denote the first infinite ordinal by $\omega$ and the first
uncountable ordinal by $\omega_1$. Throughout the paper, we use ordinal as opposed to cardinal arithmetic
in such expressions as $\omega^2$ and $\omega^\omega$. The well-known ordinal $\omegaCK$, named for Church
and Kleene, is the supremum of the recursive ordinals (those that are the order type of a recursive
relation on $\omega$). The ordinal $\omegaCK$ is also the least {\it admissible} ordinal, meaning that the
$\omegaCK$ level of G\"odel's constructible universe $L_\omegaCK$ satisfies the Kripke-Platek (KP) axioms
of set theory. We denote by $\omega_1^x$ the supremum of the $x$-recursive ordinals, and this is the same
as the least $x$-admissible ordinal, meaning that $L_{\omega_1^x}[x]\satisfies KP$. An ordinal $\alpha$ is
{\df clockable} if there is a computation of the form $\varphi_p(0)$ taking exactly $\alpha$ many steps to
halt (meaning that the $\alpha^\th$ step moves into the {\it halt} state). A {\df writable} real is one
that is the output of a computation $\varphi_p(0)$. An ordinal is writable when it is coded by a writable
real. The supremum of the writable ordinals is denoted $\lambda$, and by \cite{Welch2000:LengthsOfITTM}
this is equal to the supremum of the clockable ordinals. A real is {\df accidentally writable} when it
appears on one of the tapes at same stage during a computation of the form $\varphi_p(0)$. The supremum of
the accidentally writable ordinals, those that are coded by an accidentally writable real, is denoted
$\Sigma$. A real is {\df eventually writable} if there is a computation of the form $\varphi_p(0)$ such
that beyond some ordinal stage the real is written on the output tape (the computation need not halt).
Ordinals coded by such reals are also said to be eventually writable, and we denote the supremum of the
eventually writable ordinals by $\zeta$. Results in \cite{HamkinsLewis2000:InfiniteTimeTM} establish that
$\lambda<\zeta<\Sigma$ and that $\lambda$ and $\zeta$ are admissible. Welch \cite{Welch2000:LengthsOfITTM}
established that every computation $\varphi_p(0)$ either halts before $\lambda$ or else repeats the $\zeta$
configuration at $\Sigma$, in a transfinitely repeating loop. Furthermore, these ordinals are optimal in
the sense that the universal computation that simulates all $\varphi_p(0)$ simultaneously first enters its
repeating loop at $\zeta$, first repeating it at $\Sigma$. It follows that $\Sigma$ is not admissible.

The research in this article was initiated by the first author in a preliminary paper, which was
subsequently refined and expanded into the current three-author collaboration.

\Section Defining the Complexity Classes

Let us quickly recall the definitions of the complexity classes.

Schindler \cite{Schindler:Pnot=NP} generalized the class of polynomial decidable sets to the infinite time
context with the natural observation that every input $x\in\R$ has length $\omega$, and so the sets in $P$
should be those that are decidable in fewer steps than a polynomial function of $\omega$. Since all such
polynomials are bounded by those of the form $\omega^n$ for $n\in\omega$, he defined for $A\of\R$ that
$A\in P$\/ when there is a infinite time Turing machine $T$ and a natural number $n$ such that $T$ decides
$A$ and $T$ halts on every input in fewer than $\omega^n$ many steps.

The corresponding nondeterministic class was defined by $A\in\NP$\/ if there is an infinite time Turing
machine $T$ and a natural number $n$ such that $x\in A$ if and only if there is $y\in\R$ such that $T$
accepts $(x,y)$, and $T$ halts on every input in fewer than $\omega^n$ many steps. Sets in $\NP$\/ are
therefore simply the projections of sets in $P$.

The class $P$\/ occupies a floor a little ways upwards in the skyscraper hierarchy of classes $P_\alpha$,
indexed by the ordinals, where $A\in P_\alpha$ if and only if there is a Turing machine $T$ and an ordinal
$\beta<\alpha$ such that $T$ decides $A$, and $T$ halts on every input in fewer than $\beta$ many steps. In
this notation, the polynomial class $P$\/ is simply $P_{\omega^\omega}$, while the hierarchy continues up
through the countable ordinals to $P_{\omega_1}$, the class of sets that are decidable uniformly by some
countable stage, and $P_{\omega_1+1}$, the class of all decidable sets. We admit that the term
``polynomial'' and the letter $P$ are perhaps only appropriate at the level of $P_{\omega^\omega}$, as one
might naturally view $P_{\omega^2}$ instead as the ``linear time'' sets, $P_{\omega^{\omega^2}}$ as the
``exponential time'' sets, $P_{\epsilon_0}$ as the ``super-exponential time'' sets, $P_{\omegaCK}$ as the
``computable time'' sets, and so on, though at some point (probably already well exceeded) such analogies
become strained. Nevertheless, we retain the symbol $P$\/ in $P_\alpha$ as suggesting the polynomial time
context of classical complexity theory, because we have placed limitations on the lengths of allowed
computations. After all, infinite time Turing machines can profitably use computations of any countable
length, and so any uniform restriction to a particular countable $\alpha$ is a severe limitation. Since all
these classes concern infinite computations, one should not regard them as feasible in any practical sense.

One defines the nondeterministic hierarchy in a similar manner: $A\in\NP_\alpha$ if there is a Turing
machine $T$ and $\beta<\alpha$ such that $x\in A$ if and only if there is $y\in\R$ such that $T$ accepts
$(x,y)$, and $T$ halts on every input in fewer than $\beta$ steps. In this notation, $\NP$\/ is
$\NP_{\omega^\omega}$. Clearly, the sets in $\NP_\alpha$ are simply the projections of sets in $P_\alpha$.

As is usual in the classical context, for the nondeterministic classes we assume that the witness $y$ is
provided on a separate input tape, rather than coded together with $x$ on one input tape. This is necessary
because in order to know $P_\alpha\of\NP_\alpha\intersect\coNP_\alpha$, one wants to be able to ignore the
witness $y$ without needing extra steps of computation. When $\alpha$ is a limit ordinal or the successor
of a limit ordinal, however, one can easily manage without an extra input tape, because there is plenty of
time to decode both $x$ and the verifying witness $y$ from one input tape.

So far, these complexity classes treat every input equally in that they impose uniform bounds on the
lengths of computation, independently of the input. But it may seem more natural to allow a more
complicated input to have a longer computation. For this reason, taking $\omega_1^x$ as a natural measure
of the complexity of $x$, Schindler defined $A\in P^+$ when there is an infinite time Turing machine
deciding $A$ and halting on input $x$ in fewer than $\omega_1^x$ many steps. The corresponding
nondeterministic class is defined by $A\in\NP^+$ when there is an infinite time Turing machine $T$ such
that $x\in A$ if and only if there is $y\in\R$ such that $T$ accepts $(x,y)$, and $T$ halts on input
$(x,y)$ in fewer than $\omega_1^x$ many steps. Because this bound depends only on $x$ and not on $y$, one
can't conclude immediately that $\NP^+$ is the projection of $P^+$. One of the surprising results of the
analysis, however, is that the apparent extra power of allowing computations on input $x$ to go up to
$\omega_1^x$, as opposed to merely $\omegaCK$, actually provides no advantage (see the discussion following
Theorem \ref{NP+=Sigma^1_1}). Consequently, $\NP^+$ is the projection of $P^+$ after all.

Allowing computations to proceed a bit longer, Schindler defined that $A\in P^{++}$ when there is an
infinite time Turing machine deciding $A$ and halting on input $x$ in at most $\omega_1^x+\omega$ many
steps. Similarly, $A\in\NP^{++}$ when there is an infinite time Turing machine $T$ such that $x\in A$ if
and only if there is $y\in\R$ such that $T$ accepts $(x,y)$, and $T$ halts on any input $(x,y)$ in at most
$\omega_1^x+\omega$ many steps.

Finally, Schindler observed that any function $f$ from $\R$ to the ordinals can be viewed as bounding a
complexity class, namely, $A\in P^f$ if there is an infinite time Turing machine deciding each $x\in A$ in
fewer than $f(x)$ many steps.\footnote{This definition differs from that in \cite{HamkinsWelch2003:Pf=NPf},
which allows $\leq f(x)$ many steps in order to avoid the inevitable $+1$ that occurs when defining such
classes as $P^+$ and $P^{++}$. Here we use the original definition of \cite{Schindler:Pnot=NP}, which is
capable of describing more classes.} And $A\in\NP^f$ when there is an infinite time Turing machine $T$ such
that $x\in A$ if and only if there is $y\in\R$ such that $T$ accepts $(x,y)$, and $T$ halts on any input
$(x,y)$ in fewer than $f(x)$ many steps. In this notation, $P^+$ is the class $P^{f_0}$, where
$f_0(x)=\omega_1^x+1$, and $P^{++}=P^{f_1}$, where $f_1(x)=\omega_1^x+\omega+1$.

\Section Proving $P\not=\NP\intersect \coNP$

We begin with the basic result separating $P$\/ from $\NP\intersect\coNP$. In subsequent results we will
improve on this and precisely characterize the set $\NP\intersect\coNP$.

\Theorem. $P\not=\NP\intersect \coNP$\/ for infinite time Turing machines.\label{PisProperInNPcoNP}

\Proof: Clearly $P$ is contained in $\NP$\/ and closed under complements, so it follows that $P\of
\NP\intersect \coNP$. We now show that the inclusion is proper. Consider the halting problem for
computations halting before $\omega^\omega$ given by
$$h_{\omega^\omega}=\set{p\st \varphi_p(p)\hbox{ halts in fewer than $\omega^\omega$ steps}}.$$
We claim that $h_{\omega^\omega}\notin P$. This follows from \cite[Theorem
4.4]{HamkinsLewis2000:InfiniteTimeTM} and is an instance of Lemma \ref{h_alphaNotInPalpha} later in this
article, but let us quickly give the argument. If we could decide $h_{\omega^\omega}$ in time before
$\omega^\omega$, then we could compute the function $f(p)=1$, if $p\notin h_{\omega^\omega}$, diverge
otherwise, and furthermore we could compute this function in time before $\omega^\omega$ for input $p\notin
h_{\omega^\omega}$. If this algorithm for computing $f$ is carried out by program $q$, then $q\notin
h_{\omega^\omega}$ if and only if $f(q)\converges=1$, which holds if and only if $\varphi_q(q)$ halts in
fewer than $\omega^\omega$ steps, which holds if and only if $q\in h_{\omega^\omega}$, a contradiction.

Let us now show that $h_{\omega^\omega}\in \NP$. The idea of the proof is that the question of whether
$p\in h_{\omega^\omega}$ can be verified by inspecting (a code for) the computation sequence of
$\varphi_p(p)$ up to $\omega^\omega$. Specifically, to set this up, fix a recursive relation $\trianglelt$
on $\omega$ having order type $\omega^\omega$ and a canonical computable method of coding infinite
sequences of reals as reals, so that we may interpret any real $z$ as an infinite sequence of reals
$\<z_n\st n\in\omega>$. By combining this coding with the relation $\trianglelt$, we may view the index $n$
as representing the ordinal $\alpha$ of its order type with respect to $\trianglelt$, and we have a way to
view any real $z$ as an $\omega^\omega$-sequence of reals $\<(z)_\alpha\st\alpha<\omega^\omega>$. This
coding is computable in the sense that given any $n\in\omega$ representing $\alpha$ with respect to
$\trianglelt$, we can uniformly compute any digit of $(z)_\alpha$.

Now consider the algorithm accepting input $(p,z)$ exactly when with respect to the above coding the real
$z$ codes a halting sequence of snapshots $\<(z)_\alpha\st\alpha<\omega^\omega>$ of the computation
$\varphi_p(p)$. That is, first, each $(z)_\alpha$ codes the complete configuration of an infinite time
Turing machine, including the contents of the tapes, the position of the head, the state and the program;
second, the snapshot $(z)_{\alpha+1}$ is computed correctly from the previous snapshot $(z)_\alpha$, taking
the convention that the snapshots should simply repeat after a halt; third, the limit snapshots
$(z)_\lambda$ for limit ordinals $\lambda$ are updated correctly from the previous snapshots $(z)_\alpha$
for $\alpha<\lambda$; and finally, fourth, one of the snapshots shows the computation to have halted. Since
all of these requirements form ultimately merely an arithmetic condition on the code $z$, they can be
checked by an infinite time Turing machine in time uniformly before $\omega^2$. And since $p\in
h_{\omega^\omega}$ if and only if the computation sequence for $\varphi_p(p)$ halts before $\omega^\omega$,
we conclude that $p\in h_{\omega^\omega}$ exactly if there is a real $z$ such that $(p,z)$ is accepted by
this algorithm. Thus, $h_{\omega^\omega}\in \NP$.

To see that $h_{\omega^\omega}\in \coNP$, we simply change the fourth requirement to check that none of the
snapshots show the computation to have halted. This change means that the input $(p,z)$ will be accepted
exactly when $z$ codes a sequence of snapshots of the computation $\varphi_p(p)$, exhibiting it not to have
halted in $\omega^\omega$ many steps. Since there is a real $z$ like this if and only if $p\notin
h_{\omega^\omega}$, it follows that the complement of $h_{\omega^\omega}$ is in $\NP$, and so
$h_{\omega^\omega}\in \coNP$.\QED

Because the verification algorithm needed only to check an arithmetic condition, the argument actually
establishes $h_{\omega^\omega}\in \NP_{\omega^2}\intersect \coNP_{\omega^2}$. A closer analysis reveals
that the requirements that need to be checked are $\Pi^0_3$ (one must check that every code for a cell at a
limit stage has the right value). And since any $\Pi^0_3$ statement can be decided in time $\omega+\omega$,
it follows that $h_{\omega^\omega}$ is in $\NP_{\omega\cdot 2+2}\intersect \coNP_{\omega\cdot 2+2}$. In
fact, a bit of thought shows that the verification idea of the proof shows that any set in $\NP$\/ can be
verified by inspecting a snapshot sequence of length $\omega^\omega$, so we may actually conclude
$\NP=\NP_{\omega\cdot 2+2}$ and $\coNP=\coNP_{\omega\cdot 2+2}$. We now push these ideas harder, down to
the (optimal) level of $\omega+2$, by asking more of our witnesses.

\Theorem. The classes $\NP_\alpha$ for $\omega+2\leq\alpha\leq\omegaCK$ are all identical to the class
$\Sigma^1_1$ of lightface analytic sets. In particular, $\NP=\NP_{\omega+2}$, and so membership in any
$\NP$\/ set can be verified in only $\omega$ many steps. Similarly, the corresponding classes
$\coNP_\alpha$ are all identical to the $\Pi^1_1$ sets. Consequently, $\NP\intersect\coNP$ is exactly the
class $\Delta^1_1$ of hyperarithmetic sets.\label{NPalpha=NP}

\Proof: The idea is to have a witness not merely of the computation sequence of a given computation, but
also of all arithmetic truths. To recognize the validity of such witnesses in $\omega$ many steps, we make
use of the following two lemmas.

\SubLemma. Any $\Pi^0_2(x)$ statement can be decided on input $x$ in $\omega$ many steps.\label{Pi02}

\Proof: To decide the truth of $\forall n\exists m\psi(n,m,x)$, where $\psi$ has only bounded integer
quantifiers, one systematically considers each $n$ in turn, searching for a witness $m$ that works with
that $n$. Each time this succeeds, move to the next $n$ and flash a master flag on and then off again. If
the flag is on at a limit, it means that infinitely many $n$ were considered, so the statement is true. If
the flag is off, it means that for some $n$ the search for a witness $m$ was never completed, so the
statement is false.\QED

\SubLemma. There is an infinite time Turing machine algorithm deciding in $\omega$ many steps on input
$(a,A)$ whether $A$ is the set of arithmetic truths in $a$.\label{ArithmeticOmega}

\Proof: It is easy to see by induction on formulas that $A\of\omega$ is the set of codes for true
arithmetic statements in $a$ (that is, using $a\of\omega$ as a predicate in the language) if and only if
the following conditions, using a recursive G\"odel coding $\gcode\psi$, are satisfied:
\begin{enumerate}
\item If $\psi$ is atomic, then $\gcode\psi\in A$ if and only if $\psi$ is true.
\item $\gcode{\neg\psi}\in A$ if and only if $\gcode\psi\notin A$.
\item $\gcode{\psi\wedge\phi}\in A$ if and only if $\gcode\psi\in A$ and $\gcode\phi\in A$.
\item $\gcode{\exists u\psi(u)}\in A$ if and only if there is a natural number $n$ such that
$\gcode{\psi(n)}\in A$.
\end{enumerate}
The first three of these conditions are primitive recursive in $(a,A)$, while the fourth has complexity
$\Pi^0_2$ in $(a,A)$, making the overall complexity $\Pi^0_2$ in $(a,A)$. It follows from Lemma \ref{Pi02}
that whether or not $(a,A)$ satisfies these four conditions can be checked in $\omega$ many steps. More
concretely, we can describe an algorithm: we systematically check that $A$ satisfies each of the conditions
by considering each G\"odel code in turn. For a fixed formula, the first three conditions can be checked in
finite time. For the fourth condition, given a code for $\psi(n)$ in $A$, the algorithm can check whether
the code for $\exists u\psi(u)$ is in $A$; conversely, given that $\exists u\psi(u)$ is in $A$, let the
algorithm search for an $n$ such that $\psi(n)$ is in $A$. The point, as in Lemma \ref{Pi02}, is that if
this search fails, then at the limit one can reject the input without more ado, since it has failed
Condition (iv). Otherwise, a witness $n$ is found in finitely many steps, and the next formula is
considered.\QED

Returning to the proof of Theorem~\ref{NPalpha=NP}, we now prove that when
$\omega+2\leq\alpha\leq\omegaCK$, the classes $\NP_\alpha$ are identical. Since this is clearly a
nondecreasing sequence of classes, it suffices to show $\NP_\omegaCK\of \NP_{\omega+2}$. For this, consider
any set $B\in\NP_\omegaCK$. By definition, this means that there is a program $p$ and a recursive ordinal
$\beta$ such that $\varphi_p(x,y)$ halts in time $\beta$ for all input and $x\in B$ if and only if there is
a $y$ such that $\varphi_p$ accepts $(x,y)$. Fix a recursive relation on $\omega$ having order type
$\beta$. Consider the algorithm that accepts input $(x,y,z,A)$ exactly when $A$ codes the set of arithmetic
truths in $(x,y,z)$ and $z$ codes the computation sequence of $\varphi_p(x,y)$ of length $\beta$ (using the
fixed recursive relation for $\beta$ as the underlying order of the snapshots), and this computation
sequence shows the computation to have accepted the input. We claim that this algorithm halts in just
$\omega$ many steps.  To see this, observe first that the latter part of the condition, about $z$ coding
the computation sequence for $\varphi_p(x,y)$, is arithmetic in $(x,y,z)$. Therefore, by trusting
momentarily that $A$ is correct, it can be verified in finitely many steps by simply checking whether the
G\"odel code of that arithmetic condition is in $A$. After this, one can verify in $\omega$ many steps that
$A$ is in fact correct using the algorithm of Lemma \ref{ArithmeticOmega}. So altogether we can decide
whether $(x,y,z,A)$ has these properties in just $\omega$ many steps. And since $x\in B$ if and only if
$\varphi_p$ accepts $(x,y)$, and this happens just in case $(x,y,z,A)$ is accepted by our algorithm, where
$z$ codes the computation sequence of length $\beta$ and $A$ codes the arithmetic truths in $(x,y,z)$, we
conclude that $B\in \NP_{\omega+2}$, as desired. We have therefore proved
$\NP_{\omegaCK}\of\NP_{\omega+2}$, and so the classes $\NP_\alpha$ are identical for
$\omega+2\leq\alpha\leq\omegaCK$.

We now draw the remaining conclusions stated in the theorem. Since $\NP$\/ simply denotes
$\NP_{\omega^\omega}$, falling right in the middle of the range, it follows that $\NP=\NP_{\omega+2}$, and
so membership in any $\NP$\/ set can be verified in $\omega$ many steps. By \cite[Theorem 2.7]
{HamkinsLewis2000:InfiniteTimeTM} we know that $P_\omegaCK=\Delta^1_1$. It follows immediately that
$\NP_\omegaCK=\Sigma^1_1$, as these sets are the projections of sets in $P_\omegaCK$. So
$\NP_\alpha=\Sigma^1_1$ whenever $\omega+2\leq\alpha\leq\omegaCK$, as these classes are all identical. And
finally, by taking complements, we conclude as well that $\coNP_\alpha=\Pi^1_1$ whenever
$\omega+2\leq\alpha\leq\omegaCK$.\QED

It will follow from Theorem \ref{NP+=Sigma^1_1} that this result can be extended at least one more step, to
$\omegaCK+1$, because $\NP_{\omegaCK}=\NP_{\omegaCK+1}$.

\Corollary. $\NP\not=\coNP$\/ for infinite time Turing machines.

\Proof: The classes $\Sigma^1_1$ and $\Pi^1_1$ are not identical.\QED

Both Theorems \ref{PisProperInNPcoNP} and \ref{NPalpha=NP} can also be proved using the model-checking
technique of \cite{HamkinsWelch2003:Pf=NPf}, which we will use extensively later in this article.

\Section Proving $P^+=\NP^+\intersect\coNP^+$

At first glance, the class $P^+$ appears much more generous than the earlier classes, because computations
on input $x$ are now allowed up to $\omega_1^x$ many steps, which can be considerably larger than
$\omegaCK$. But it will follow from Theorem \ref{NP+=Sigma^1_1} that if a set is in $P^+$, then there is an
algorithm deciding it in uniformly fewer than $\omegaCK$ many steps, much sooner than required. Our
arguments rely on the following fact from descriptive set theory.

\Lemma. $\Pi^1_1$ absorbs existential quantification over $\Delta^1_1$. That is, if $B$ is $\Pi^1_1$ and $x
\in A\iff \exists y \in \Delta^1_1(x)\, B(x,y)$, then $A$ is $\Pi^1_1$ as
well.\label{AbsorbingExistentialQ}

\Proof: This lemma is a special case of \cite[Theorem 4D.3]{Moschovakis1980:DescriptiveSetTheory}, and is
due to Kleene. We provide a proof sketch here. Let $U$ be a universal $\Pi^1_1$ set and suppose
$y\in\Delta^1_1(x)$. Then there is an integer $i_0$ such that $y(n)=m$ if and only if $U(i_0,x,n,m)$. Let
$U^*$ be a $\Pi^1_1$ set uniformizing $U$, so that for all $i,x,n$ if there is an $m$ with $U(i,x,n,m)$
then there is a unique $m$ with $U^*(i,x,n,m)$. In particular, $y(n)=m$ if and only if $U^*(i_0,x,n,m)$. So
we have altogether that $x \in A$ if and only if there is an integer $i$ such that for all $n$ there is
exactly one $m$ with $U^*(i,x,n,m)$, and for all $y$, either $B(x,y)$ or there are $n,m$ with
$U^*(i,x,n,m)$ and $y(n) \not= m$. In other words, we say that there is an index $i$ of a computation of a
real $z$ via a $\Pi^1_1(x)$ recursive function such that $B(x,z)$. As $\Pi^1_1$ is closed under
quantification over integers, this shows that $A$ is in $\Pi^1_1$, as desired.\QED

\Theorem. \
\begin{enumerate}
\item $\NP^+ = \Sigma^1_1=\NP=\NP_\alpha$ whenever $\omega+2\leq\alpha\leq\omegaCK+1$.
\item $P^+=\Delta^1_1=P_{\omegaCK}=P_{\omegaCK+1}$.
\item $P^+ = \NP^+\intersect \coNP^+$.
\end{enumerate}
\label{NP+=Sigma^1_1}

\Proof: For (i), we have already proved in Theorem~\ref{NPalpha=NP} that $\Sigma^1_1=\NP$, and since
clearly $\NP\of \NP^+$, it follows that $\Sigma^1_1\of \NP^+$. Conversely, suppose that $A\in \NP^+$. This
means that there is an infinite time Turing machine program $p$ such that $\varphi_p(x,y)$ halts on all
input $(x,y)$ in fewer than $\omega_1^x$ many steps, and $x\in A$ if and only if there is a real $y$ such
that $\varphi_p$ accepts $(x,y)$. The set $A$ is therefore the projection of the set $$B=\set{(x,y)\st
\varphi_p{\rm\ accepts\ }(x,y)}.$$ In order to see that $A$ is in $\Sigma^1_1$, it suffices to show
$B\in\Sigma^1_1$ (and our argument shows just as easily that $B\in\Delta^1_1$). The complement of $B$ is
the set $\neg B=\set{(x,y)\st \varphi_p{\rm\ rejects\ }(x,y)}$, and these computations also have length
less than $\omega_1^x$. It follows that the computation sequence for $\varphi_p(x,y)$ exists in the model
$L_{\omega_1^x}[x,y]$, and so $(x,y)\in\neg B$ if and only if $L_{\omega_1^x}[x,y]\satisfies\theta(x,y)$,
where $\theta(x,y)$ asserts that the computation $\varphi_p(x,y)$ rejects the input. Since this is a
$\Sigma_1$ assertion, it follows that $(x,y)\in\neg B$ if and only if there is an ordinal
$\beta<\omega_1^x$ such that $L_\beta[x,y]\satisfies\theta(x,y)$. Since the model $L_\beta[x,y]$ is
hyperarithmetic in $(x,y)$, and any well-founded model showing the computation to reject the input will do,
we see that $(x,y)\in\neg B$ if and only if there is a real $z\in\Delta^1_1(x,y)$ coding a well-founded
model of $V=L[x,y]$ that satisfies $\theta(x,y)$. Since the property of coding a well-founded model (of any
theory) is $\Pi^1_1$ in the theory, it follows by Lemma \ref{AbsorbingExistentialQ} that $\neg B$ is
$\Pi^1_1$. Consequently, $B\in\Sigma^1_1$, and so $A$, being the projection of $B$, is in $\Sigma^1_1$ as
well. So we have proved that $\NP^+=\Sigma^1_1$. It follows from Theorem \ref{NPalpha=NP} that
$\NP^+=\NP_\alpha$ whenever $\omega+2\leq\alpha\leq\omegaCK$. The remaining case of $\alpha=\omegaCK+1$
follows from (ii) and the observation that $\NP_\omegaCK=\NP_{\omegaCK+1}$, as these are the projections of
$P_\omegaCK=P_{\omegaCK+1}$.

For (ii) and (iii), observe that since $\NP^+=\Sigma^1_1$, it follows that $\coNP^+=\Pi^1_1$, and so
$P^+\of \NP^+\intersect \coNP^+=\Sigma^1_1\intersect \Pi^1_1=\Delta^1_1$, which by \cite[Theorem
2.7]{HamkinsLewis2000:InfiniteTimeTM} is equal to $P_{\omegaCK}$, which is a subset of $P_{\omegaCK+1}$,
which is clearly a subset of $P^+$. So all of them are equal, as we claimed.\QED

The fact that $P^+=\Delta^1_1$ was Theorem 2.13 of \cite{Schindler:Pnot=NP}, and one can view our argument
here as a detailed expansion of that argument. In fact, however, once one knows
$P^+=P_{\omegaCK}=\Delta^1_1$, it follows immediately that sets in $\NP^+$ are projections of sets in
$P_\omegaCK=\Delta^1_1$, since all the computations halt uniformly before $\omegaCK$, which is certainly
not larger than $\omega_1^x$, and consequently $\NP^+=\Sigma^1_1$. By this means, Theorem
\ref{NP+=Sigma^1_1} follows directly from \cite[Theorem 2.13]{Schindler:Pnot=NP}.

The fact that $P^+=P_\omegaCK$ should be surprising---and we mentioned this in the introduction---because
it means that although the computations deciding $x\in A$ for $A\in P^+$ are allowed to compute up to
$\omega_1^x$, in fact there is an algorithm needing uniformly fewer than $\omegaCK$ many steps. So the
difference between $\omegaCK$ and $\omega_1^x$, which can be substantial, gives no advantage at all in
computation. An affirmative answer to the following question would explain this phenomenon completely.

\Question. Suppose an algorithm halts on each input $x$ in fewer than $\omega_1^x$ steps. Then does it halt
uniformly before $\omegaCK$?

Secondly, the fact that $P_\omegaCK=P_{\omegaCK+1}$ is itself surprising, because the difference in the
definitions of these two classes is exactly the difference between requiring the computations to halt
before $\omegaCK$ and requiring them to halt {\it uniformly} before $\omegaCK$, that is, before some fixed
$\beta<\omegaCK$ on all input. Since the classes $P_\omegaCK=P_{\omegaCK+1}$ are equal, any set that can be
decided before $\omegaCK$ can be decided uniformly before $\omegaCK$.

Finally, let us close this section with a more abstract view of Theorem \ref{NP+=Sigma^1_1}. Suppose that
$f:\R\to\omega_1$ is Turing invariant and for some $\Sigma_1$ formula $\varphi$ we have $f(x)=\alpha$ if
and only if $L[x]\satisfies\varphi(x,\alpha)$. We define the pointclass $\Gamma^f$ by $A\in\Gamma^f$ if and
only if there is some $\Sigma_1$ formula $\theta$ such that $x \in A\iff L_{f(x)}[x] \satisfies \theta(x)$.
Then for ``natural" f one should be able to show that $P^f = \Delta^f = \NP^f \intersect \coNP^f$ and
$\NP^f = \Gamma^f$-dual, by our arguments above. (We do not attempt to classify the functions $f$ for which
these equations hold true.) For $f(x) = \omega_1^x$, these equations collapse to the the statement of
Theorem \ref{NP+=Sigma^1_1}. The pointclasses $\Gamma^f$ exhaust all of $\Delta^1_2$.



\Section The Question Whether $P_\alpha=\NP_\alpha\intersect\coNP_\alpha$

We turn now to the relation between $P_\alpha$ and $\NP_\alpha\intersect \coNP_\alpha$ for various ordinals
$\alpha$. We begin with the observation that the classes $P_\alpha$ increase with every clockable limit
ordinal $\alpha$.

\Definition. The lightface {\df halting problem} is the set $h=\{\,p\st \varphi_p(p)$ halts$\}$.
Approximating this, for any ordinal $\alpha$ the {\df halting problem for $\alpha$} is the set
$h_\alpha=\{\,p\st \varphi_p(p)$ halts in fewer than $\alpha$ many steps$\}$. We sometimes denote
$h_{\alpha+1}$ by $h_\leqalpha$, to emphasize the fact that it is concerned with computations of length
less than or equal to $\alpha$. Similarly, we denote $P_{\alpha+2}$, $\NP_{\alpha+2}$ and
$\coNP_{\alpha+2}$ by $P_\leqalpha$, $\NP_\leqalpha$ and $\coNP_\leqalpha$, respectively, as these classes
are also concerned only with the computations of length less than or equal to $\alpha$. It follows that
$\NP_\leqalpha$ is the projection of $P_\leqalpha$, and $\coNP_\leqalpha$ consists of the complements of
sets in $\NP_\leqalpha$.

\Lemma. If $\alpha$ is any ordinal, then $h_\alpha\notin P_\alpha$. Indeed, $h_\alpha\notin P_{\alpha+1}$.
In particular, $h_{\leqalpha}\notin P_\leqalpha$. However, if $\alpha$ is a clockable limit ordinal, then
$h_\alpha\in P_\leqalpha$.\label{h_alphaNotInPalpha}

\Proof: Suppose to the contrary that $h_\alpha\in P_{\alpha+1}$ for some ordinal $\alpha$. It follows that
there is an algorithm $q$ deciding $h_\alpha$ in fewer than $\alpha$ many steps. That is, the computation
of $\varphi_q(x)$ halts in fewer than $\alpha$ many steps on any input and $x\in h_\alpha$ if and only if
$\varphi_q$ accepts $x$. Consider the modified algorithm $q_0$ that runs $\varphi_q(x)$, but when the
algorithm is just about to move into the {\it accept} state, it instead jumps into a non-halting
transfinite repeating loop. This algorithm computes a function $\varphi_{q_0}(x)$ which halts in fewer than
$\alpha$ steps if $x\notin h_\alpha$ and diverges otherwise. Therefore, $q_0\in h_\alpha$ if and only if
$\varphi_{q_0}(q_0)$ halts, which holds if and only if $q_0\notin h_\alpha$, a contradiction. So we have
established $h_\alpha\notin P_{\alpha+1}$ for any ordinal $\alpha$. It follows, in particular, that
$h_\leqalpha=h_{\alpha+1}\notin P_{\alpha+2}=P_\leqalpha$.

Finally, when $\alpha$ is a clockable limit ordinal, consider the algorithm that on input $p$ simulates
both the computation $\varphi_p(p)$ and the $\alpha$ clock (simulating $\omega$ many steps of each in every
$\omega$ many actual steps). If the computation stops before the clock runs out, the algorithm accepts the
input, but if the clock runs out, it rejects the input. By placing the first column of the clock's
computation in the actual first column, the algorithm will be able to detect that the clock has stopped at
exactly stage $\alpha$, and thereby halt in $\alpha$ steps. So $h_\alpha\in P_\leqalpha$.\QED

\Corollary. If $\alpha$ is a clockable limit ordinal, then $P_\alpha\ofnoteq
P_\leqalpha$.\label{PalphaIncreases}

For recursive ordinals $\alpha$ and even $\alpha\leq\omegaCK+1$, the question whether
$P_\alpha=\NP_\alpha\intersect\coNP_\alpha$ is already settled by Theorem~\ref{NPalpha=NP}, and we
summarize the situation here.

\Theorem. $P_\alpha\not=\NP_\alpha\intersect \coNP_\alpha$ whenever $\omega+2\leq\alpha<\omegaCK$. Equality
is attained at $\omegaCK$ and its successor with $$P_\omegaCK=\NP_\omegaCK\intersect
\coNP_\omegaCK=\Delta^1_1=P_{\omegaCK+1}=\NP_{\omegaCK+1}\intersect\coNP_{\omegaCK+1}.$$\label{Palpha}

\vskip-20pt\Proof: For $\alpha<\omegaCK$ we know by Corollary \ref{PalphaIncreases} that $P_\alpha$ is a
proper subset of $P_\omegaCK$, which by Theorem \ref{NPalpha=NP} is equal to
$\NP_\omegaCK\intersect\coNP_\omegaCK=\NP_\alpha\intersect\coNP_\alpha$. So none of the earlier classes
$P_\alpha$ are equal to $\NP_\alpha\intersect \coNP_\alpha$; but at the top we do achieve the equalities
$P_\omegaCK=\NP_\omegaCK\intersect \coNP_\omegaCK$ and
$P_{\omegaCK+1}=\NP_{\omegaCK+1}\intersect\coNP_{\omegaCK+1}$ because by Theorem \ref{NP+=Sigma^1_1} these
are both instances of the identity $\Delta^1_1=\Sigma^1_1\intersect\Pi^1_1$.\QED

Corollary \ref{PalphaIncreases} and Theorem \ref{Palpha} show that the class $\Delta^1_1$ of
hyperarithmetic sets is ramified by the increasing hierarchy $\union_{\alpha<\omegaCK}P_\alpha$ in a way
similar to the traditional hyperarithmetic hierarchy $\Delta^1_1=\union_{\alpha<\omegaCK}\Delta^0_\alpha$,
and one can probably give a tight analysis of the interaction of these two hierarchies.

We now prove that the pattern of Theorems \ref{NPalpha=NP} and \ref{Palpha}---where the classes
$\NP_\alpha$ are identical for $\alpha$ in the range from $\omega+2$ up to $\omegaCK+1$---is mirrored
higher up, within any contiguous block of clockable ordinals. It will follow that $P_\alpha$ is properly
contained in $\NP_\alpha\intersect\coNP_\alpha$ within any such block of clockable ordinals. We
subsequently continue the pattern at the top of any such block, by proving that
$P_\beta=\NP_\beta\intersect\coNP_\beta$ for the ordinal $\beta$ that begins the next gap in the clockable
ordinals.

\Theorem. If $[\nu,\beta)$ is a contiguous block of infinite clockable ordinals, then all the classes
$\NP_\alpha$ for $\nu+2\leq\alpha\leq\beta+1$ are identical. Consequently, all the corresponding classes
$\coNP_\alpha$ for such $\alpha$ are identical as well.\label{NPalphaIdenticalInContiguousBlock}

\Proof: Since the sequence of classes $\NP_\alpha$ is nondecreasing, it suffices to show
$\NP_{\beta+1}\of\NP_{\nu+2}$. Suppose $B\in\NP_{\beta+1}$, so that there is an algorithm $e$ such that
$\varphi_e(x,y)$ halts on every input in in time less than $\beta$, and $x\in B$ if and only if there is
$y$ such that $\varphi_e$ accepts $(x,y)$. Since $\nu$ is clockable, there is a program $q_0$ such that
$\varphi_{q_0}(0)$ takes exactly $\nu$ steps to halt.

Consider the algorithm which on input $(x,z)$ checks, first, whether $z$ codes a model $M_z$ of $KP$
containing $x$ in which the computation $\varphi_{q_0}(0)$ halts and there is a $y\in M_z$ such that
$\varphi_e$ accepts $(x,y)$; and second, verifies that $M_z$ is well-founded up to $\nu^{M_z}$, the length
of the clock computation $\varphi_{q_0}(0)$ in $M_z$. If both of these requirements are satisfied, then the
algorithm accepts the input, and otherwise rejects it.

If $x\in B$, then there is a $y$ such that $\varphi_e$ accepts $(x,y)$, and so we may choose $z$ coding a
fully well-founded model $M_z$ that is tall enough to see this computation and $\varphi_{q_0}(0)$. It
follows that $(x,z)$ will be accepted by our algorithm. Conversely, if $(x,z)$ is accepted by our
algorithm, then the corresponding model $M_z$ is well-founded up to $\nu$. Since the well-founded part is
admissible and no clockable ordinal is admissible, $M_z$ must be well-founded beyond the length of the
computation $\varphi_e(x,y)$ (which is less than $\beta$), since all the ordinals in $[\nu,\beta)$ are
clockable. Therefore, $M_z$ will have the correct (accepting) computation for $\varphi_e(x,y)$, and so
$x\in B$. Thus, our algorithm nondeterministically decides $B$. And as before, since $\nu$ is inadmissible,
this algorithm will either discover ill-foundedness below $\nu$, halting in time at most $\nu$, or else
halt at $\nu$ with well-foundedness up to $\nu^{M_z}=\nu$. So $B\in \NP_{\nu+2}$, as desired.

The corresponding fact for $\coNP_\alpha$ follows by taking complements.\QED

\Corollary. $P_\alpha\not=\NP_\alpha\intersect\coNP_\alpha$ for any clockable ordinal $\alpha\geq\omega+2$,
except possibly when $\alpha$ ends a gap in the clockable ordinals or is the successor of such a gap-ending
ordinal.\label{PalphaForClockableAlpha}

\Proof: If $\alpha\geq\omega+2$ is clockable but is neither a gap-ending ordinal nor the successor of a
gap-ending ordinal, then there is an infinite ordinal $\nu<\nu+2\leq\alpha$ such that $[\nu,\alpha]$ is a
contiguous block of clockable ordinals. By Theorem \ref{NPalphaIdenticalInContiguousBlock}, the classes
$\NP_\xi$ are identical for $\nu+2\leq\xi\leq\beta+1$, where $\beta$ is the next admissible beyond
$\alpha$. Since by Corollary \ref{PalphaIncreases} the corresponding classes $P_\xi$ increase at every
clockable limit ordinal in this range and are subsets of $\NP_\xi\intersect\coNP_\xi$, it follows that
$P_\alpha\ofnoteq \NP_\alpha\intersect\coNP_\alpha$.\QED

\Corollary. In particular, $P_\leqalpha\not=\NP_\leqalpha\intersect\coNP_\leqalpha$ for any infinite
clockable ordinal $\alpha$.\label{PleqalphaForClockableAlpha}

\Proof: This is an instance of the previous theorem, because $P_\leqalpha=P_{\alpha+2}$ and $\alpha+2$ is
neither a limit ordinal nor the successor of a limit ordinal.\QED

Because of the possible exceptions in Corollary \ref{PalphaForClockableAlpha} at the gap-ending ordinals,
we do not have a complete answer to the following question.

\Question. Is $P_\alpha\not=\NP_\alpha\intersect\coNP_\alpha$ for any clockable ordinal
$\alpha\geq\omega+2$?

The first unknown instances of this occur at the first gap-ending ordinal $\omegaCK+\omega$ and its
successor $\omegaCK+\omega+1$. Thus, we don't know whether
$$P_{\omegaCK+\omega}=\NP_{\omegaCK+\omega}\intersect\coNP_{\omegaCK+\omega},$$ nor do we know whether
$$P_{\omegaCK+\omega+1}=\NP_{\omegaCK+\omega+1}\intersect\coNP_{\omegaCK+\omega+1}.$$ A related question
concerns the gap-starting ordinals and their successors, such as $\omegaCK$ and $\omegaCK+1$, where we have
proved the equalities
$$P_\omegaCK=\NP_\omegaCK\intersect\coNP_\omegaCK\qquad{\rm and}\qquad
P_{\omegaCK+1}=\NP_{\omegaCK+1}\intersect\coNP_{\omegaCK+1}.$$
 We will now show that this phenomenon is completely general, appealing to the following unpublished
results of Philip Welch.

\Lemma.({\cite[Lemma 2.5]{Welch:ActionOfOneTapeMachines}}) If $\alpha$ is a clockable ordinal, then every
ordinal up to the next admissible beyond $\alpha$ is writable in time
$\alpha+\omega$.\label{QuicklyWritable}

\Theorem.({\cite[Theorem 1.8]{Welch:ActionOfOneTapeMachines}}) Every ordinal beginning a gap in the
clockable ordinals is admissible.\label{GapStartersAreAdmissible}

This latter result is a converse of sorts to \cite[Theorem 8.8]{HamkinsLewis2000:InfiniteTimeTM}, which
establishes that no admissible ordinal is clockable. It is not the case, however, that the gap-starting
ordinals are exactly the admissible ordinals below $\lambda$, because admissible ordinals can appear in the
middle of a gap. To see that this phenomenon occurs, observe that the suprema of the writable and
eventually writable ordinals are both admissible, with no clockable ordinals in between, and this situation
reflects downwards into an actual gap, because an algorithm can search for accidentally writable admissible
ordinals having no clockable ordinals in between, and halt when they are found.

\Theorem. Suppose that $\beta$ begins a gap in the clockable ordinals. Then
$P_\beta=\NP_\beta\intersect\coNP_\beta$. Furthermore, if $\beta$ is in addition not a limit of
non-clockable ordinals, then
$P_\beta=P_{\beta+1}=\NP_{\beta+1}\intersect\coNP_{\beta+1}$.\label{Gap-starters}

\Proof: Let us suppose first that $\beta$ begins a gap in the clockable ordinals, but is not a limit of
non-clockable ordinals, so that there is some $\nu<\beta$ such that $[\nu,\beta)$ is a contiguous block of
clockable ordinals. Since $\nu$ is clockable, it follows by Lemma \ref{QuicklyWritable} that there is a
real $u$ coding $\nu$ that is writable in time $\nu+\omega$, which is of course still less than $\beta$. We
claim that $\beta=\omega_1^u$. To see this, observe that since $\beta$ is admissible, $L_\beta$ has the
computation producing $u$ and so $u\in L_\beta$. Consequently, $\beta$ is $u$-admissible and so
$\omega_1^u\leq\beta$. Conversely, since $u$ codes $\nu$ and there are no admissible ordinals in
$[\nu,\beta)$, it follows that $\beta\leq\omega_1^u$, and so $\beta=\omega_1^u$.

Next, we relativize Theorem \ref{Palpha} with respect to an oracle for $u$, concluding that
$P_{\omega_1^u}^u=\NP_{\omega_1^u}^u\intersect\coNP_{\omega_1^u}^u=\Delta^1_1(u)=
P_{\omega_1^u+1}^u=\NP_{\omega_1^u+1}^u\intersect\coNP_{\omega_1^u+1}^u$, where the superscript indicates
the presence of an oracle for $u$. But since $u$ is writable in time $\nu+\omega<\beta$ by Lemma
\ref{QuicklyWritable}, we can simulate such an oracle by simply taking the time first to write it out. By
admissibility, $\nu+\omega+\beta=\beta$, and so this preparatory step will not cause any ultimate delay in
our calculations. Therefore, $P_\beta=P_\beta^u$, $\NP_\beta=\NP_\beta^u$ and $\coNP_\beta=\coNP_\beta^u$,
and the same for $\beta+1$. We conclude that
$P_\beta=\NP_\beta\intersect\coNP_\beta=\Delta^1_1(u)=P_{\beta+1}=\NP_{\beta+1}\intersect\coNP_{\beta+1}$,
as desired.

It remains to consider the case of gap-starting ordinals $\beta$ that are limits of gaps. In this case,
$\beta$ is a limit of ordinals $\xi$ that begin gaps but are not limits of non-clockable ordinals (they
begin the ``successor'' gaps), and consequently by the previous paragraph satisfy
$P_\xi=\NP_\xi\intersect\coNP_\xi=P_{\xi+1}=\NP_{\xi+1}\intersect\coNP_{\xi+1}$. Because $P_\beta$ is the
union of the nondecreasing sequence of classes $P_\xi$ for $\xi<\beta$, and the same for $\NP_\beta$ and
$\coNP_\beta$, it follows that
$$P_\beta=\Union_{\xi<\beta}P_\xi=\Union_{\xi<\beta}\NP_\xi\intersect\coNP_\xi=
(\Union_{\xi<\beta}\NP_\xi)\intersect(\Union_{\xi<\beta}\coNP_\xi)=\NP_\beta\intersect\coNP_\beta,$$ and so
the proof is complete.\QED

\Corollary. In particular, $P_\lambda=\NP_\lambda\intersect\coNP_\lambda$, where $\lambda$ is the supremum
of the clockable ordinals.

More generally, we ask for a characterization of these exceptional ordinals.

\Question. Exactly which ordinals $\alpha$ satisfy $P_\alpha=\NP_\alpha\intersect\coNP_\alpha$?

Just for the record, let us settle the question for the classes $P_\omega$ and $P_{\omega+1}$, as well as
$P_n$ for finite $n$, which are all trivial in the sense that they involve only finite computations. The
class $P_\omega$ concerns the uniformly finite computations, while $P_{\omega+1}$ allows arbitrarily long
but finite computations. The class $P_n$ for finite $n$ concerns computations having at most $n-2$ steps.
Observe that $P_0=P_1=\emptyset$ because computations have nonnegative length, and
$P_2=\{\,\R,\emptyset\,\}$ because a computation halts in $0$ steps only when the {\it start} state is
identical with either the {\it accept} or {\it reject} states. Infinite computations first appear with the
class $P_{\omega+2}$.

\Theorem. For the classes corresponding to finite computations:\
\begin{enumerate}
 \item $P_n=\NP_n=\coNP_n$ for any finite $n$. Consequently, $P_n=\NP_n\intersect\coNP_n$.
 \item $P_\omega=\NP_\omega=\coNP_\omega$. Consequently, $P_\omega=\NP_\omega\intersect\coNP_\omega$.
 \item $P_{\omega+1}=\Delta^0_1$, $\NP_{\omega+1}=\Sigma^0_1$ and
$\coNP_{\omega+1}=\Pi^0_1$. Consequently, $P_{\omega+1}=\NP_{\omega+1}\intersect \coNP_{\omega+1}$.
\end{enumerate}\label{Pomega}

\Proof: For (i), observe that the computations putting a set in $P_n$, $\NP_n$ or $\coNP_n$ are allowed at
most $n-2$ many steps, and so the sets they decide must depend on at most the first $n-2$ digits of the
input. But any such set is in $P_n$, because if membership in $A\of\R$ depends on the first $n-2$ digits of
the input, then there is a program which simply reads those digits, remembering them with states, and moves
to the {\it accept} or {\it reject} states accordingly. So $P_n=\NP_n=\coNP_n$. Claim (ii) follows, because
$P_\omega=\Union_nP_n=\Union_n\NP_n=\NP_\omega$.

For (iii), observe that a set $B$ is in $P_{\omega+1}$ if $x\in B$ can be decided by a Turing machine
program that halts in finitely many steps. Since this is precisely the classical notion of (finite time)
computability, it follows that $P_{\omega+1}=\Delta^0_1$, the recursive sets of reals. If $B\in
\NP_{\omega+1}$, there is an algorithm $p$ such that $\varphi_p(x,y)$ halts in finitely many steps on all
input and $x\in B$ if and only if there is a $y$ such that $\varphi_p$ accepts $(x,y)$. Thus, $x\in B$ if
and only if there is a finite piece $y\restrict n$ such that $\varphi_p$ accepts $(x,y\restrict n)$, where
the piece is long enough that the algorithm never inspects $y$ beyond $n$ bits. Since this has now become
an existential quantifier over the integers, we conclude that $B\in \Sigma^0_1$. Conversely, every set in
$\Sigma^0_1$ is clearly the projection of a set in $\Delta^0_1$, so we conclude
$\NP_{\omega+1}=\Sigma^0_1$. By taking complements, $\coNP_{\omega+1}=\Pi^0_1$.\QED

Returning our focus to the infinite computations, let us now consider the case of ordinals that are not
necessarily clockable. Our first observation is that the key idea of the proof of Theorem
\ref{PisProperInNPcoNP}---the fact that one could easily recognize codes for $\omega^\omega$ or any other
recursive ordinal---generalizes to the situation where one has only nondeterministic algorithms for
recognizing the ordinals in question.

\Definition. An ordinal $\alpha$ is {\df recognizable} (in time $\xi$) when there is a nonempty set of
reals coding $\alpha$ that is decidable (in time $\xi$). The ordinal $\alpha$ is {\df nondeterministically
recognizable} (in time $\xi$) if there is a nonempty set of codes for $\alpha$ that is nondeterministically
decidable (in time $\xi$).

If $\alpha$ is nondeterministically recognizable in time $\xi$, then the set $\WO_\alpha$ of {\it all}
reals coding $\alpha$ is nondeterministically decidable in time $\xi$, because a real is in $\WO_\alpha$ if
and only if there is an isomorphism from the relation it codes to the relation coded by any other real
coding $\alpha$.

\Lemma. If an ordinal $\alpha$ is nondeterministically recognizable in time $\xi$, then
$h_\alpha\in\NP_\leqxi\intersect\coNP_\leqxi$.\label{Recognizable}

\Proof: Suppose that $\alpha$ is nondeterministically recognizable in time $\xi$, so there is a nonempty
set $D$ of codes for $\alpha$ that is in $\NP_\leqxi$. We may assume both $\alpha$ and $\xi$ are at least
$\omegaCK$. Consider the algorithm that on input $(p,u,v,w)$ checks, first, that $u$ codes a linearly
ordered relation on $\omega$ with respect to which $v$ codes the snapshot sequence of $\varphi_p(p)$,
showing it to halt, and second, that $(u,w)$ is accepted by the nondeterministic algorithm deciding $D$,
verifying $u\in D$. If $p\in h_\alpha$, then the computation $\varphi_p(p)$ halts in fewer than $\alpha$
many steps, and so we may choose a real $u\in D$ coding $\alpha$, along with a real $w$ witnessing that
$u\in D$, and a real $v$ coding the halting snapshot sequence of $\varphi_p(p)$, so that $(p,u,v,w)$ is
accepted by our algorithm. Conversely, if $(p,u,v,w)$ is accepted by our algorithm, then because $(u,w)$
was accepted by the algorithm for $D$, we know $u$ really codes $\alpha$, and so the snapshot sequence must
be correct in showing $\varphi_p(p)$ to halt before $\alpha$, so $p\in h_\alpha$. Finally, the algorithm
takes $\xi$ steps, because the initial check takes fewer than $\omega^2$ steps, being arithmetic, and so
the computation takes $\omega^2+\xi=\xi$ many steps altogether. Thus, $h_\alpha\in\NP_\xi$.

To see that $h_\alpha\in\coNP_\xi$, simply modify the algorithm to check that $v$ codes a snapshot sequence
with respect to the relation coded by $u$, but $v$ shows the computation {\it not} to halt.\QED

One can use the same idea to show that if $\NP_\alpha$ contains a set of codes for ordinals unbounded in
$\alpha$, then $P_\alpha\not= \NP_\alpha$.

We will now apply this result to show that $P_\alpha\not=\NP_\alpha\intersect\coNP_\alpha$ for all
sufficiently large countable ordinals $\alpha$. Recall from the introduction that $\lambda<\zeta<\Sigma$
refer to the suprema of the writable, eventually writable and accidentally writable ordinals, respectively.
The first two of these are admissible, while the latter is not, and every computation either halts before
$\lambda$ or repeats the $\zeta$ configuration at $\Sigma$. And furthermore, $\Sigma$ is characterized by
being the first repeat point of the universal computation simulating all $\varphi_p(0)$ simultaneously.

\Theorem. If $\Sigma+2\leq\alpha$, then $P_\alpha\not=\NP_\alpha\intersect\coNP_\alpha$. In fact, the class
$\NP_{\leq\Sigma}\intersect\coNP_{\leq\Sigma}$ contains a nondecidable set, the halting problem
$h$.\label{PleqSigma}

\Proof: The proof relies on the following.

\SubLemma. $\Sigma$ is nondeterministically recognizable in time $\Sigma$.\label{SigmaIsNDRecognizable}

\Proof: The model-checking algorithm of \cite[Theorem 1.7]{HamkinsWelch2003:Pf=NPf} essentially shows this,
but let us sketch the details here. By results in \cite{Welch2000:LengthsOfITTM}, the ordinal $\Sigma$ is
the first stage at which the universal computation (simulating $\varphi_p(0)$ for all programs $p$) repeats
itself. Consider the algorithm which on input $(x,y)$ checks whether $x$ codes a relation on $\omega$ and
$y$ codes a model $M_y\satisfies``KP+\Sigma$ exists'' containing $x$ and satisfying the assertion that the
order type of $x$ is $\Sigma^{M_y}$. If $y$ passes this test, then the algorithm counts-through the
relation coded by $x$ to verify that it is well-founded. If all these tests are passed, then the algorithm
accepts in the input, and otherwise rejects it. If the well-founded part of $M_y$ exceeds the true
$\Sigma$, then $M_y$ will have the correct value for $\Sigma$, and the algorithm will take exactly $\Sigma$
many steps. If the well-founded part of $M_y$ lies below $\Sigma$, then this will be discovered before
$\Sigma$ and the algorithm will halt before $\Sigma$. Finally, because $\Sigma$ is not admissible, the
well-founded part of $M_y$ cannot be exactly $\Sigma$, and so in every case our algorithm will halt in at
most $\Sigma$ many steps. And since the acceptable $x$ have order type $\Sigma$, this shows that
$\WO_\Sigma$, the set of reals coding $\Sigma$, is nondeterministically decidable in $\Sigma$ steps, as
desired.\QED

By Lemma \ref{Recognizable}, it follows that $h_\Sigma\in\NP_\leqSigma\intersect\coNP_\leqSigma$. But since
$\Sigma$ is larger than every clockable ordinal, it follows that $h_\Sigma=h$, the full lightface halting
problem. So we have established that if $\Sigma+2\leq\alpha$, then the halting problem $h$ is in
$\NP_\alpha\intersect\coNP_\alpha$. Since $h$ is not decidable, it cannot be in $P_\alpha$. So
$P_\alpha\not=\NP_\alpha\intersect\coNP_\alpha$.\QED

This establishes $P_\alpha\not=\NP_\alpha\intersect\coNP_\alpha$ for all but countably many $\alpha$. We
close this section with another definition and an application.

\Definition. An ordinal $\alpha$ is {\df nondeterministically clockable} if there is an algorithm $p$ which
halts in time at most $\alpha$ for all input and in time exactly $\alpha$ for some input. More generally,
$\alpha$ is {\df nondeterministically clockable before $\beta$} if there is an algorithm that halts before
$\beta$ on all input and in time exactly $\alpha$ for some input.

Such an algorithm can be used as a clock for $\alpha$ in nondeterministic computations, since there are
verifying witnesses making the clock run for exactly the right amount of time, with a guarantee that no
other witnesses will make the clock run on too long.

\Theorem. If $\alpha$ is an infinite nondeterministically clockable limit ordinal, then
$P_\leqalpha\not=\NP_\leqalpha$.\label{PleqalphaForNClockableAlpha}

\Proof: By Lemma \ref{h_alphaNotInPalpha}, it follows that $h_{\alpha+\omega}\notin P_\leqalpha$. But we
claim that $h_{\alpha+\omega}\in\NP_\leqalpha$. By Theorem \ref{NPalpha=NP} we may assume
$\alpha>\omegaCK$, because when $\alpha$ is recursive $h_{\alpha+\omega}$ is hyperarithmetic, and hence
already in $\NP_{\omega+2}$. Fix a nondeterministic clock for $\alpha$, a program $e$ such that
$\varphi_e(z)$ halts in exactly $\alpha$ many steps for some $z$ and in at most $\alpha$ many steps on all
other input. We will now nondeterministically decide $h_{\alpha+\omega}$ by the following algorithm. On
input $(x,y,z)$, first determine whether $x$ is some finite $p$. If not, then reject the input, otherwise,
check whether $y$ codes a model $M_y$ of $KP$ containing $z$ and satisfying the assertion that
$\varphi_e(z)$ halts, with $\varphi_p(p)$ halting at most finitely many steps later. Since this is an
arithmetic condition on $y$, it can be checked in fewer than $\omega^2$ many steps. Next, assuming that
these tests have been passed successfully, we verify that the model $M_y$ is well-founded up to what it
thinks is the halting time of $\varphi_e(z)$, which we denote $\alpha^{M_y}$. If ill-foundedness is
discovered, we reject the input. By flashing a master flag every time we delete what is the current
smallest (in the natural ordering of $\omega$) element still in the field, we can tell at a limit stage
that we have finished counting, and when this occurs, we accept the input.

Let's argue that this algorithm accomplishes what we want. First of all, if $p\in h_{\alpha+\omega}$, then
$\varphi_p(p)$ halts before $\alpha+\omega$ and there is a real $z$ such that $\varphi_e(z)$ halts in
$\alpha$ steps and a real $y$ coding a fully well-founded model $M_y\satisfies KP$ in which these
computations exist. So the previous algorithm will accept the input $(p,y,z)$. Conversely, if the algorithm
accepts $(p,y,z)$ for some $y$ and $z$, then the corresponding model $M_y$ is well-founded up to the length
of the computation $\varphi_e(z)$, which is at most $\alpha$ because the computation $\varphi_e(z)$ in
$M_y$ agrees with the actual computation as long as the model remains well-founded. It follows that model
is also well-founded for an additional $\omega$ many steps, and so the model is correct about
$\varphi_p(p)$ halting before $\alpha+\omega$. So the algorithm does nondeterministically decide
$h_{\alpha+\omega}$.

It remains to see that the algorithm halts in at most $\alpha$ many steps on all input. Since
$\omegaCK\leq\alpha$, it follows that $\omega^2+\alpha=\alpha$, and so the initial checks of those
arithmetic properties do not ultimately cause any delay. The only question is how many steps it takes to
check the well-foundedness of $M_y$ up to $\alpha^{M_y}$. If $M_y$ is well-founded up to $\alpha^{M_y}$,
then this takes exactly $\alpha^{M_y}$ many steps (as the count-through algorithm is designed precisely to
take $\beta$ steps to count through a relation of limit order type $\beta$), and this is at most $\alpha$.
If $M_y$ is ill-founded below $\alpha^{M_y}$, then this will be discovered exactly $\omega$ many steps
beyond the well-founded part of $M_y$, and so the algorithm will halt in at most $\alpha$ many steps.
Lastly, the well-founded part of $M_y$ cannot be exactly $\alpha$, because $\alpha$ is not $z$-admissible.
So in any case, on any input the algorithm halts in at most $\alpha$ many steps.\QED

This argument does not seem to establish that $P_\leqalpha\not=\NP_\leqalpha\intersect\coNP_\leqalpha$ for
such $\alpha$, however, because one cannot seem to use a nondeterministic clock in this algorithm to verify
that a computation $\varphi_p(p)$ has {\it not} halted. The problem is that a prematurely halting
nondeterministic clock might cause the algorithm to think that $\varphi_p(p)$ does not halt in time
$\alpha+\omega$ even when it does, which would lead to false acceptances for the complement of
$h_{\alpha+\omega}$.

\Section The Cases of $P^f$ and $P^{++}$

Let us turn now to the question of whether $P^f=\NP^f\intersect\coNP^f$, where $f:\R\to\ORD$. A special
case of this is the question of whether $P^{++}=\NP^{++}\intersect\coNP^{++}$, because $P^{++}=P^{f_1}$,
where $f(x)=\omega_1^x+\omega+1$. We consider only functions $f$ that are {\df suitable}, meaning that
$f(x)\leq_T f(y)$ whenever $x\leq_T y$ and $f(x)\geq\omega+1$.

Many of the instances of the question whether $P^f=\NP^f\intersect\coNP^f$ are actually solved by a close
inspection of the arguments of \cite{HamkinsWelch2003:Pf=NPf}, though the results there were stated only as
$P^f\not=\NP^f$. The point is that the model-checking technique of verification used in those arguments is
able to verify both positive and negative answers.

But more than this, the next theorem shows that the analysis of whether $P^f=\NP^f\intersect\coNP^f$, at
least for sets of natural numbers, reduces to the question of whether
$P_\alpha=\NP_\alpha\intersect\coNP_\alpha$, where $\alpha=f(0)+1$. And since the previous section provides
answers to this latter question for many values of $\alpha$, we will be able to provide answers to the
former question as well, in Corollaries \ref{f(0)>Sigma} and \ref{f(q)isClockable}.

\Theorem. For any suitable function $f$ and any set $A$ of natural numbers,
\begin{enumerate}
 \item $A\in P^f$ if and only if $A\in P_{f(0)+1}$;
 \item $A\in\NP^f$ if and only if $A\in\NP_{f(0)+1}$;
 \item $A\in\coNP^f$ if and only if $A\in\coNP_{f(0)+1}$.
\end{enumerate}\label{P_f(0)}

\Proof: By suitability, $f(0)\leq f(x)$ for all $x$, and $f(0)=f(n)$ for all natural numbers $n$. Since any
set in $P_{f(0)+1}$ is decided by an algorithm that takes fewer than $f(0)$ many steps, it follows that
$P_{f(0)+1}\of P^f$. Conversely, suppose that $A\of\omega$ and $A\in P^f$. So there is an algorithm that
decides whether $x\in A$ in fewer than $f(x)$ many steps. Although this algorithm might be allowed to take
many steps on a complicated input $x$ for which $f(x)$ may be large, we know since $A\of\omega$ that the
ultimate answer will be negative unless $x\in\omega$. Thus, we design a more efficient algorithm by
rejecting any input $x$ that does not code a natural number. Since the natural number $n$ is coded by the
sequence consisting of a block of $n$ ones, followed by zeros, the sequences that don't code natural
numbers are precisely the sequence of all ones, plus those containing the substring $01$. While continuing
with the algorithm to decide $A$, our modified algorithm searches for the substring $01$ in the input, and
also turns on a flag if $0$ is encountered in the input. This algorithm decides $n\in A$ in fewer than
$f(n)=f(0)$ many steps, and rejects all other input either in finitely many steps, if the input contains
$01$, or in $\omega$ many steps, if the input has no zeros. It therefore places $A$ in $P_{f(0)+1}$, as
desired.

A similar argument establishes the result for $\NP^f$ and $\NP_{f(0)+1}$. Specifically, if $A\in\NP^f$,
then there is a nondeterministic algorithm such that $x\in A$ if and only if the algorithm accepts $(x,y)$
for some $y$. Once again, we can modify this algorithm to reject any input $(x,y)$ in finitely many steps
unless $x$ codes some finite $n$, in which case the algorithm is carried out as before. The result is that
$x\in A$ is decided in finite time unless $x=n\in\omega$, in which case it is decided in fewer than
$f(n)=f(0)$ many steps, placing $A$ in $\NP_{f(0)+1}$. The result for $\coNP^f$ and $\coNP_{f(0)+1}$
follows by taking complements.\QED

The argument of \cite[Theorem 3.1]{HamkinsWelch2003:Pf=NPf} essentially proves the following result, though
that result is stated merely as $P^f\not=\NP^f$. Here, we will derive it as a corollary to the previous
theorem and Theorem~\ref{PleqSigma}. Note that if $f:\R\to\ORD$ is suitable, then $f(q)=f(0)$ for any
finite $q$.

\Corollary. If $f:\R\to\ORD$ is suitable and $f(0)>\Sigma$, then $P^f$ is properly contained in
$\NP^f\intersect\coNP^f$.\label{f(0)>Sigma}

\Proof: This follows immediately from Theorems \ref{PleqSigma} and \ref{P_f(0)}, because the halting
problem $h$, being a set of natural numbers and in $\NP_{{\leq}\Sigma}\intersect\coNP_{{\leq}\Sigma}$, must
be in $\NP^f\intersect\coNP^f$, but it is not decidable and consequently not in $P^f$.\QED

\Corollary. If $f:\R\to\ORD$ is suitable and $f(0)$ is clockable, but does not end a gap in the clockable
ordinals, then $P^f$ is properly contained in $\NP^f\intersect\coNP^f$.\label{f(q)isClockable}

\Proof: By Theorem \ref{P_f(0)}, the sets of natural numbers in $P^f$ and $\NP^f\intersect\coNP^f$ are
exactly those in $P_{\alpha+1}$ and $\NP_{\alpha+1}\intersect\coNP_{\alpha+1}$, respectively, where
$\alpha=f(0)$. Since $f(0)$ does not end a gap in the clockable ordinals, it follows that $\alpha+1$ is
neither a gap-ending ordinal nor the successor of a gap-ending ordinal. Therefore, by Corollary
\ref{PalphaForClockableAlpha} there are sets of natural numbers in
$\NP_{\alpha+1}\intersect\coNP_{\alpha+1}$ that are not in $P_{\alpha+1}$. Consequently, there are sets of
natural numbers in $\NP^f\intersect\coNP^f$ that are not in $P^f$.\QED

An instance of this settles the question for $P^{++}$.

\Corollary. $P^{++}\not=\NP^{++}\intersect\coNP^{++}$.

\Proof: This follows from Corollary \ref{f(q)isClockable} and the fact that $P^{++}=P^{f_1}$, where
$f_1(x)=\omega_1^x+\omega+1$. By \cite[Theorem 3.2]{HamkinsLewis2000:InfiniteTimeTM}, the ordinal
$\omegaCK+\omega$ is clockable, and consequently so is $\omegaCK+\omega+1$.\QED

So the previous corollaries establish that $P^f\not=\NP^f\intersect\coNP^f$ for many or most functions $f$.
But of course, we have examples of ordinals $\alpha$ for which $P_\alpha=\NP_\alpha\intersect\coNP_\alpha$,
such as $\alpha=\omegaCK$ or $\alpha=\omegaCK+1$. If $f$ is the constant function $f(x)=\omegaCK$, then it
is easy to see that  $P^f=P_{\omegaCK+1}$ and $\NP^f=\NP_{\omegaCK+1}$, and this provides an example where
$P^f=\NP^f\intersect\coNP^f$, even when $P^f\not=\NP^f$. The equation $P^+=\NP^+\intersect\coNP^+$ provides
another such example.

\bigskip\bigskip
\noindent
Vinay Deolalikar\\
Hewlett-Packard Research\\
1501 Page Mill Road, M/S 3U-4, Palo Alto, CA 94304\\
Tel: (650) 857 8605, Fax: (650) 852 3791\\
http://www.hpl.hp.com/personal/Vinay\!\_Deolalikar/

\bigskip
\noindent
Joel David Hamkins\\
Georgia State University {\it\&} The City University of New York\\
The College of Staten Island of CUNY and The CUNY Graduate Center\\
Mathematics Program, 365 Fifth Avenue, New York, NY 10016\\
http://jdh.hamkins.org

\bigskip
\noindent
Ralf-Dieter Schindler\\
Institut f\"ur formale Logik, Universit\"at Wien\\
1090 Wien, Austria\\
rds@logic.univie.ac.at\\
http://www.logic.univie.ac.at/${}^\sim$rds/

\bibliographystyle{alpha}
\bibliography{MathBiblio,HamkinsBiblio}

\end{document}